\documentclass[11,openbib]{article}
\usepackage[T1]{fontenc}
\usepackage{amsmath}
\usepackage{amsfonts}
\usepackage{amssymb}

\newcommand{\comp}{\raisebox{0pt}{$\scriptstyle \circ \, $}}

\allowdisplaybreaks

\begin{document}
\begin{center}
{\Large \textbf{Aronsza{\Large j}n and Sikorski subcartesian \\ 
\vspace{.01in} differential spaces}}\\[0pt]
\mbox{}\vspace{-0.2in} \\[0pt]
\mbox{}\\[0pt]
Richard Cushman and J\k{e}drzej \'{S}niatycki\footnotemark \\[0pt]
\mbox{}\\[0pt]
\end{center}
\footnotetext{Department of Mathematics and Statistics, University of Calgary \newline
e-mail: r.h.cushman@gmail.com and sniatycki@gmail.com} 
\hfill \today

\begin{abstract}
We construct a natural transformation between the category of Aronszajn subcartesian 
spaces and the category of subcartesian differential spaces, which is a subcategory of 
Sikorski differential spaces. 
\end{abstract} \bigskip

Differential spaces were introduced in 1967 by Sikorski \cite{sikorski1}. A
comprehensive presentation of Sikorski's theory of differential spaces is
contained in his book \cite{sikorski72}. For the current state of the theory
consult \cite[chpt.1]{sniatycki}. In this theory, differential structure of
a space is given by its algebra of smooth functions. In 1967 Aronszajn \cite{aronszajn} 
introduced the notion of a subcartesian space with its smooth structure
described in terms of an atlas, analogous to the standard definition of a 
differentiable manifold. \smallskip 

We can characterize a smooth Hausdorff manifold of dimension $n$ as Hausdorff
differential space $S$ such that every point $x\in S$ has an open neighbourhood $U$ diffeomorphic to an open subset $V$ of $\mathbb{R}^{n}$. Here the
differential structures on $U$ and $V$ are generated by restrictions of
smooth functions on $S$ and $\mathbb{R}^{n}$, respectively. We can weaken this definition by not requiring that $V$ is open in $\mathbb{R}^{n}$ and allowing $n$ to be an arbitrary non-negative integer dependent on 
$U.$ This leads to the following notion of a subcartesian differential
space. A differential space $S$ is \emph{subcartesian} if it is Hausdorff
and every point $x\in S$ has a neighbourhood $U_{x}$ diffeomorphic to a
differential subspace $V_{x}$ of $\mathbb{R}^{n_{x}}$. \medskip

\noindent \textbf{Proposition 1} A differential subspace $S$ of $\mathbb{R}^{n}$ is a
subcartesian differential space. Moreover, if $h\in C^{\infty }(S)$, then 
for every point $x\in S$ there is an open neighbourhood $U$ of $x$ in 
$S\subseteq \mathbb{R}^{n}$ and a function $f\in C^{\infty }(\mathbb{R}^{n})$
such that $h_{\mid U}=f_{\mid U}$. \medskip 

\noindent \textbf{Proof.}
Since the set $S$ is a subset of $\mathbb{R}^{n}$, it is Hausdorff. Hence,
the differential subspace $S$ of $\mathbb{R}^{n}$ is subcartesian. Let $%
(x_{1},...,x_{n})$ be the natural coordinate functions on $\mathbb{R}^{n}$.
By definition of a differential structure generated by a family of functions, see 
\cite{sniatycki}, there exists $F\in C^{\infty }(\mathbb{R}^{n})$ such that
\begin{equation*}
h_{\mid U}=F(x_{1},...,x_{n})_{\mid U}.
\end{equation*}%
But $f=F(x_{1},...,x_{n})$ is in $C^{\infty }(S)$ \hfill $\square $ \medskip

A \emph{subcartesian space of Aronszajn} is a Hausdorff topological 
space $S$ endowed with an atlas 
$\mathfrak{A}=\{\varphi :U_{\varphi }\rightarrow V_{\varphi }\} $,
where $\varphi :U_{\varphi }\rightarrow V_{\varphi }$ is a homeomorphism of
an open subset $U_{\varphi }$ of $S$ onto a subset $V_{\varphi }$ of ${\mathbb{R}}^{n_{\varphi }}$, which has the following properties:\smallskip 

\indent\textbf{1}. \parbox[t]{4in}{The domains $\{U_{\varphi }\mid \varphi \in \mathfrak{A}\}$ form
an open cover of $S.$}\smallskip 

\indent \textbf{2}. \parbox[t]{4in}{For every $\varphi $, $\psi \in \mathfrak{A}$ and every $x\in
U_{\varphi }\cap U_{\psi }$, there exists a $C^{\infty }$-mapping $\textbf{s}$
extending $\psi \comp {\varphi }^{-1}: \varphi (U_{\varphi }\cap U_{\psi})\rightarrow \psi (U_{\varphi }\cap U_{\psi })$ in a neighbourhood of $\varphi (x)\in {\mathbb{R}}^{n_{\varphi }}$. Also, there exists a $C^{\infty }$-mapping $\textbf{t}$ extending $\varphi \comp {\psi }^{-1}:\psi (U_{\varphi }\cap U_{\psi })\rightarrow \varphi 
(U_{\varphi }\cap U_{\psi })$ in a neighbourhood of $\psi (x)\in {\mathbb{R}}^{n_{\psi }}$.} \medskip 

\noindent As in the theory of manifolds, there may be different atlases giving rise
to the same subcartesian structure. A subcartesian structure determines a
maximal atlas, which can be interpreted as the union of all the atlases
giving this structure.\smallskip

In \cite{walczak} Walczak characterized the largest class of differential
spaces which satisfy a condition analogous to condition 2 of the definition of Aronszajn subcartesian spaces. However, he did not give an explicit proof that subcartesian and Aronszajn subcartesian spaces are equivalent. 
The goal of this paper is to show that the definitions of subcartesian and Aronszajn subcartesian spaces are
equivalent and to construct a natural transformation between the category of Aronszajn subcartesian 
spaces and subcartesian differential spaces. \medskip

First, we give Aronszajn's definition of a smooth map between Aronszajn  
subcartesian spaces. If $(S_{1},\mathfrak{A}_{1})$ and $(S_{2},\mathfrak{A}_{2})$ are subcartesian spaces of Aronszajn, a map $\chi :S_{1}\rightarrow
S_{2}$ is \emph{smooth} if, for every $x\in S_{1}$, there exist $\varphi _{1}\in 
\mathfrak{A}_{1}$ and $\varphi _{2}\in \mathfrak{A}_{2}$ such that $x\in
U_{\varphi _{1}}$, $\chi (x)\in U_{\varphi _{2}}$ and $\varphi _{2}\comp
\chi \comp \varphi _{1}^{-1}:V_{\varphi _{1}}\rightarrow V_{\varphi _{2}}$ \linebreak
extends to a $C^{\infty }$-map of a neighbourhood of $\varphi _{1}(x)\in 
\mathbb{R}^{n_{\varphi _{1}}}$ to a neighbourhood of $\varphi _{2}(\chi
(x))\in \mathbb{R}^{n_{\varphi _{2}}}$. A map $\chi :S_{1}\rightarrow S_{2}$
is a \emph{diffeomorphism} of Aronszajn's subcartesian spaces if it is smooth,
invertible and its inverse is smooth.\smallskip 

Let $S$ be a subcartesian differential space with a differential structure $%
C^{\infty }(S)$. By definition, for every $x\in S$, there
exists a diffeomorphism $\varphi _{x}:U_{x}\subseteq S \rightarrow V_{x}\subseteq 
\mathbb{R}^{n_{x}}$. Here, $U_{x}$ is an open neighbourhood of $x$ in $S$
endowed with the differential structure $C^{\infty }(U_{x})$. Similarly, $V_{x}$
is a subset of $\mathbb{R}^{n_{x}}$ generated by pulling back the standard differential structure of 
${\mathbb{R}}^{n_x}$ by the inclusion mapping 
$V_{x}\hookrightarrow \mathbb{R}^{n_{x}}$. \medskip

\noindent \textbf{Proposition 2} Let $S$ be a differential space with $\mathfrak{A} = 
{ \{ {\varphi }_x: U_x \rightarrow V_{x} \} }_{x\in S}$. The pair $(S,\mathfrak{A})$ is an Aronszajn
subcartesian space. Moreover, the topology of $(S,\mathfrak{A})$ is the same
as the topology of the subcartesian differential space $(S,C^{\infty }(S)).$\smallskip

\noindent \textbf{Proof.} \newline 
(\textbf{a}) Let $S^{\prime \prime }$ be a differential subspace of $\mathbb{R}^{n}$. Denote by $\mathfrak{A}^{\prime \prime }$ the atlas for $S^{\prime \prime }$ defined above. In this case the inclusion map 
$\iota ^{\prime \prime }:S^{\prime \prime}\hookrightarrow \mathbb{R}^{n}$ is smooth by definition of a differential subspace. Also, $S^{\prime \prime }$ is Hausdorff topological subspace of $\mathbb{R}^{n}$. \smallskip 

For $k=i,j$ let $\varphi _{k}^{\prime \prime }:U_{k}^{\prime \prime}\rightarrow 
\varphi _{k}^{\prime \prime }(U_{k}^{\prime \prime })\in \mathbb{R}^{m_{k}}$ be two charts in 
${\mathfrak{A}}^{\prime \prime }$ such that 
$U_{ij}^{\prime \prime }=U_{i}^{\prime \prime }\cap U_{j}^{\prime \prime}\neq \varnothing $. 
Then $V_{i}^{\prime \prime }=\varphi _{i}^{\prime \prime}(U_{ij}^{\prime \prime })$ is a differential subspace 
of $\mathbb{R}^{m_{i}} $, $V_{j}^{\prime \prime }=
\varphi _{j}^{\prime \prime}(U_{ij}^{\prime \prime })$ is a differential subspace of 
${\mathbb{R}}^{m_{j}} $, and the map $\varphi _{ij}^{\prime \prime }= \varphi _{i}^{\prime
\prime }\comp (\varphi _{j}^{\prime \prime })^{-1}:V_{j}\rightarrow V_{i}$
is a diffeomorphism. Let $(z_{1},...,z_{m_{i}})$ be coordinates in $\mathbb{R}^{m_{i}}$, 
considered as maps $z_{k}:\mathbb{R}^{m_{i}}\rightarrow \mathbb{R} $, for $k=1,...,m_{i}$. 
If $\iota _{i}^{\prime \prime }:V_{i}^{\prime \prime }\rightarrow \mathbb{R}^{m_{i}}$ is the inclusion map, then, 
\begin{equation*}
\iota _{i}^{\prime \prime }\comp \varphi _{ij}^{\prime \prime }=(z_{1}\comp
\varphi _{ij}^{\prime \prime },...,z_{m_{i}}\comp \varphi _{ij}^{\prime
\prime }):V_{j}^{\prime \prime }\rightarrow \mathbb{R}^{m_{i}}
\end{equation*}%
is smooth. For $k=1,...,m_{i}$, the function $z_{k}\comp \varphi
_{ij}^{\prime \prime }\ $is in $C^{\infty }(V_{j}^{\prime \prime })$. By 
proposition 1, for $y\in V_{j}^{\prime \prime }$ there exists an open
neighbourhood $W_{k}^{\prime \prime }$ of $y$ in $V_{j}^{\prime \prime
}\subseteq \mathbb{R}^{m_{i}}$ and a function $f_{k}^{\prime \prime }\in
C^{\infty }(\mathbb{R}^{m_{i}})$ such that $(z_{k}\comp \varphi
_{ij}^{\prime \prime })_{\mid W_{k}}=f_{k\mid W_{k}}^{\prime \prime }$. Let $%
W^{\prime \prime }=W_{1}^{\prime \prime }\cap ...\cap W_{k}^{\prime \prime
}\cap ...\cap W_{m_{i}}^{\prime \prime }$. Then, $f_{k\mid W^{\prime \prime
}}^{\prime \prime }=(z_{k}\comp \varphi _{ij}^{\prime \prime })_{\mid
W^{\prime \prime }}$ for $k=1,...,m_{i}.$ Therefore,%
\begin{equation*}
(z_{1}\comp \varphi _{ij}^{\prime \prime },...,z_{m_{i}}\comp \varphi
_{ij}^{\prime \prime })_{\mid W^{\prime \prime }}=(f_{1}^{\prime \prime
},...,f_{m_{i}}^{\prime \prime })_{\mid W^{\prime \prime }}:W^{\prime \prime} \rightarrow \mathbb{R}^{m_{j}},
\end{equation*}%
and $(f_{1}^{\prime \prime },...,f_{m_{i}}^{\prime \prime }):\mathbb{R}%
^{m_{i}}\rightarrow \mathbb{R}^{m_{j}}$ is an extension of 
\begin{displaymath}
\varphi _{ij\mid
W}^{\prime \prime }=(\varphi _{i}^{\prime \prime }\comp (\varphi
_{j}^{\prime \prime })^{-1})_{\mid W}:W^{\prime \prime }\rightarrow
V_{i}^{\prime \prime } .
\end{displaymath}
Similarly, we can construct local extensions of $\varphi _{j}^{\prime \prime
}\comp (\varphi _{j}^{\prime \prime })^{-1}:V_{i}^{\prime \prime
}\rightarrow V_{j}^{\prime \prime }$. \ Therefore, for every differential
subspace $S^{\prime \prime }$ of $\mathbb{R}^{n}$, $\mathfrak{A}^{\prime \prime }$ 
is an atlas. Hence, $(S^{\prime \prime },\mathfrak{A}^{\prime \prime })$ is an Aronszajn subcartesian space. \smallskip

\par \noindent \textbf{(b)} Let $S^{\prime }$ be a differential space diffeomorphic to a
differential subspace $S^{\prime \prime }$ of $\mathbb{R}^{n}$ and let $\chi
:S^{\prime }\rightarrow S^{\prime \prime }$ be a diffeomorphism between
these spaces. $S^{\prime }$ is Hausdorff because $\chi $ is a homeomorphism. Let 
$\mathfrak{A}^{\prime \prime }=\mathfrak{\{\varphi }^{\prime \prime }%
\mathfrak{\}}$ be the atlas on $S^{\prime \prime}$ introduced in (\textbf{a}). It consists of charts $%
\varphi ^{\prime \prime }:U^{\prime \prime }\rightarrow \varphi ^{\prime
\prime }(U^{\prime \prime })\subseteq \mathbb{R}^{m}$, where $U^{\prime
\prime }$ is an open subset of $S^{\prime \prime }$ and $\varphi ^{\prime \prime }$ is a
diffeomorphism onto its image. Then 
\begin{equation*}
\mathfrak{A}^{\prime }=\{\varphi ^{\prime }\mid \varphi ^{\prime }=\chi
^{-1}\comp \varphi ^{\prime \prime }\comp \chi :U^{\prime }=\chi
^{-1}(U^{\prime \prime })\rightarrow \varphi ^{\prime }(U^{\prime
})\subseteq \mathbb{R}^{m},\mathrm{for}%
\text{ }\varphi ^{\prime \prime }\in \mathfrak{A}^{\prime \prime }\}
\end{equation*}
is an atlas on $S^{\prime }$. Since $\chi :S^{\prime }\rightarrow
S^{\prime \prime }$ is a diffeomorphism of subcartesian spaces, the atlas 
$\mathfrak{A}^{\prime }$ satisfies condition 1 of the definition of Aronszajn subcartesian space. \smallskip 

For $k=i,j$ let $\varphi _{k}^{\prime \prime }:U_{k}^{\prime \prime}\rightarrow 
\varphi _{k}^{\prime \prime }(U^{\prime \prime })_{k}\in 
\mathbb{R}^{m_{k}}$ be two charts in $\mathfrak{A}^{\prime \prime }$ such
that $U_{ij}^{\prime \prime }=U_{i}^{\prime \prime }\cap U_{j}^{\prime
\prime }\neq \varnothing $, as above. We follow the arguments given in
(\textbf{a}). Let $U_{i}^{\prime }=\chi ^{-1}(U_{i}^{\prime \prime }),$ and 
\begin{equation*}
\varphi _{i}^{\prime }=\varphi _{i}^{\prime \prime }\comp \chi
:U_{i}^{\prime }\rightarrow \varphi _{i}^{\prime }(U_{i}^{\prime })=\varphi
_{i}^{\prime \prime }\comp \chi \comp \chi ^{-1}(U_{i}^{\prime \prime
})=\varphi _{i}^{\prime \prime }(U_{i}^{\prime \prime }).
\end{equation*}
Then, 
\begin{equation*}
U_{ij}^{\prime }=U_{i}^{\prime }\cap U_{j}^{\prime }=\chi
^{-1}(U_{i}^{\prime \prime })\cap \chi ^{-1}(U_{j}^{\prime \prime })=\chi
^{-1}\left( U_{i}^{\prime \prime }\cap U_{j}^{\prime \prime }\right) =\chi
^{-1}\left( U_{ij}^{\prime \prime }\right) \neq \varnothing .
\end{equation*}%
Moreover, 
\begin{align}
\varphi _{i}^{\prime }(U_{ij}^{\prime }) & =\varphi _{i}^{\prime }\comp \chi
^{-1}\left( U_{i}^{\prime \prime }\cap U_{j}^{\prime \prime }\right)
=\varphi _{i}^{\prime \prime }\comp \chi \comp \chi ^{-1}\left(
U_{i}^{\prime \prime }\cap U_{j}^{\prime \prime }\right) \notag \\
& =\varphi _{i}^{\prime \prime }\left( U_{i}^{\prime \prime }\cap U_{j}^{\prime \prime
}\right) =\varphi _{i}^{\prime \prime }\left( U_{ij}^{\prime \prime }\right)
=V_{ij}^{\prime \prime } \notag
\end{align}%
and $\varphi _{j}^{\prime }(U_{ij}^{\prime })$ are differential subspaces of 
$\mathbb{R}^{m_{i}}$ and $\mathbb{R}^{m_{j}}$, respectively, and the map 
\begin{equation*}
\varphi _{ij}^{\prime }=\varphi _{i}^{\prime }\comp (\varphi _{j}^{\prime
})^{-1}=(\varphi _{i}^{\prime \prime }\comp \chi )\comp (\varphi
_{j}^{\prime \prime }\comp \chi )^{-1}=\varphi _{i}^{\prime \prime }\comp
(\varphi _{j}^{\prime \prime })^{-1}=\varphi _{ij}^{\prime \prime
}:V_{j}^{\prime \prime }\rightarrow V_{i}^{\prime \prime }
\end{equation*}%
is a diffeomorphism. In (\textbf{a}) we showed that this map satisfies the condition
2 of the definition of Aronszajn subcartesian space. Hence $(S^{\prime },\mathfrak{A}^{\prime }),$ where $%
S^{\prime }$ is a differential space diffeomorphic to a differential
subspace of $\mathbb{R}^{n}$, is a subcartesian space of
Aronszajn.\smallskip 

\noindent (\textbf{c}) Suppose that $S$ is a subcartesian differential space.
By definition $S$ is Hausdorff. Moreover, each point $x\in S$
has an open neigbhourhood $U_{x}\subseteq S,$ and a diffeomorphism $\varphi
_{x}:U_{x}\rightarrow V_{x}$, where $V_{x}$ is a differential subspace of $%
\mathbb{R}^{n_{x}}$, for some $n_{x}\in \mathbb{Z}_{\geq 0}.$  By the
arguments in (\textbf{b}), each differential subspace $U_{x}$ of $S$ is an Aronszajn subcartesian space. 
Since $S$ is Hausdorff, it suffices to show that the maps ${\varphi }_x$ also satisfy the second condition of
the definition of Aronszajn subcartesian space. However, this condition is local, and is satisfied in each
open neighbourhood $U_{x}$ of the covering $\{U_{x}\}$ of $S$. Hence, the
Hausdorff subcartesian differential space $S$ with the atlas $\mathfrak{A}%
=\{\varphi _{x}:U_{x}\rightarrow V_{x}\}_{x\in S}$ is a subcartesian space
of Aronszajn.\smallskip\ 

The topology of the subcartesian differential space $S$ is encoded in its
differential structure. By definition, $S$ is Hausdorff. Moreover, the domains of charts in $\mathfrak{A}$ form an open cover of $S$. Hence the the topology on the Aronszajn
subcartesian space $(S,\mathfrak{A})$ coincides with the topology of the subcartesian 
differential space $(S,C^{\infty }(S)).$ \hfill $\square $ \smallskip 

\vspace{-.15in}For $i=1,2,$ let $S_{i}$ be a subcartesian differential space with
differential structure $C^{\infty }(S_{i})$, and let $\mathfrak{A}_{i}$ be the set of diffeomorphisms of 
open subsets of $S_{i}$ onto subsets of Euclidean spaces. By proposition 2, $(S_{1},\mathfrak{A}_{1})$
and $(S_{2},\mathfrak{A}_{2})$ are subcartesian spaces of Aronszajn. \medskip 

\noindent \textbf{Proposition 3} If $\chi :(S_{1},C^{\infty
}(S_{1}))\rightarrow (S_{2},C^{\infty }(S_{2}))$ is smooth map of
differential spaces and the atlasses $\mathfrak{A}_{1}$ and $\mathfrak{A}_{2} $ are are constructed in terms of the differential structures $C^{\infty }(S_{1})$ and $C^{\infty }(S_{1})$ as in proposition 2, then  
$\chi :(S_{1},\mathfrak{A}_{1})\rightarrow (S_{2},\mathfrak{A}_{2})$ is a smooth map of Aronszajn subcartesian
spaces.\medskip

\noindent \textbf{Proof.}
Suppose that $\chi :S_{1}\rightarrow S_{2}$ is a smooth map of differential
spaces. This means that $\chi ^{\ast }f\in C^{\infty }(S_{1})$ for
each $f\in C^{\infty }(S_{2}).$ We want to show that $\chi
:S_{1}\rightarrow S_{2}$ is smooth in the sense of Aronszajn subcartesian differential
space. According to the definition of Aronszjan subcartesian space, we have to show that for every $x\in
S_{1}$, there exist $\varphi _{1}:U_{\varphi _{1}}\rightarrow V_{\varphi
_{1}}$ in $\mathfrak{A}_{1}$ and $\varphi _{2}:U_{\varphi _{2}}\rightarrow
V_{\varphi _{2}}$ in $\mathfrak{A}_{2}$, such that $x\in U_{\varphi _{1}}$, $%
\chi (x)\in U_{\varphi _{2}}$ and $\varphi _{2}\comp \chi \comp \varphi
_{1}^{-1}:V_{\varphi _{1}}\rightarrow V_{\varphi _{2}}$ extends to a $%
C^{\infty }$-map of a neighbourhood of $\varphi _{1}(x)\in \mathbb{R}%
^{n_{\varphi _{1}}}$ to a neighbourhood of $\varphi _{2}(\chi (x))\in 
\mathbb{R}^{n_{\varphi _{2}}}$. In order to simplify the notation, we write $%
n_{1}=n_{\varphi _{1}}$ and $n_{2}=n_{\varphi _{2}}$.\smallskip

For each $i=1,2,$ the domains of charts in $\mathfrak{A}_{i}$ cover $S_{i}$.
Hence, for $x\in S_{1}$, there exists $\varphi _{1}\in \mathfrak{A}_{1}$
such that $x\in U_{\varphi _{1}},$ and there exists $\varphi _{2}\in 
\mathfrak{A}_{2}$ such that $\chi (x)\in U_{\varphi _{2}}.$ For $i=1,2,$ $%
\varphi _{i}:U_{\varphi _{i}}\rightarrow V_{\varphi _{i}}$ is a
diffeomorphism of a differential subspace $U_{\varphi _{i}}$ of $S_{i}$ onto
a differential subspace $V_{\varphi _{i}}$ of $\mathbb{R}^{n_{i}}$. Since $%
\chi :S_{1}\rightarrow S_{2}$ is smooth in the sense of differential spaces,
it follows that $\varphi _{2}\comp \chi \comp \varphi _{1}^{-1}:V_{\varphi
_{1}}\rightarrow V_{\varphi _{2}}$ is a smooth map of a differential
subspace $V_{\varphi _{1}}$ of $\mathbb{R}^{n_{1}}$ to a differential
subspace $V_{\varphi _{2}}$ of $\mathbb{R}^{n_{2}}$, and $\varphi _{2}\comp
\chi \comp \varphi _{1}^{-1}(\varphi _{1}(x))=\varphi _{2}(\chi (x)).$%
\smallskip

As in the proof of proposition 2, let $(z_{1},...,z_{n_{2}})$ be
coordinates in $\mathbb{R}^{n_{2}}$, considered as maps $z_{k}:\mathbb{R}^{n_{2}}\rightarrow \mathbb{R}$, for $k=1,\ldots ,n_{2}$. If ${\iota}_{2}:V_{\varphi _{2}}\rightarrow \mathbb{R}^{n_{2}}$ is the inclusion map,
then, 
\begin{equation*}
{\iota }_{2}\comp \varphi _{2}\comp \chi \comp \varphi _{1}^{-1}= 
(z_{1}\comp \varphi _{2}\comp \chi \comp \varphi _{1}^{-1},...,z_{n_{2}}\comp {\varphi}_{2}\comp \chi \comp 
{\varphi }_{1}^{-1}):V_{\varphi _{1}}\rightarrow \mathbb{R}^{n_{2}}
\end{equation*}%
is smooth. For $k=1,...,n_{2}$, the function $z_{k}\comp \varphi _{2}\comp
\chi \comp \varphi _{1}^{-1}$ lies in $C^{\infty }(V_{\varphi _{1}})$. By
proposition 1, for $y\in V_{\varphi _{1}},$ there exists an open
neighbourhood $W_{k}$ of $y$ in $V_{\varphi _{1}}\subseteq \mathbb{R}^{n_{1}} $ and a function 
$f_{k}\in C^{\infty }(\mathbb{R}^{n_{1}})$ such
that $(z_{k}\comp \varphi _{2}\comp \chi \comp \varphi _{1}^{-1})_{\mid
W_{k}}=f_{k\mid W_{k}}$. Let $W=W_{1}\cap ...\cap W_{k}\cap ...\cap
W_{n_{2}} $. Then, $f_{k\mid W}=(z_{k}\comp \varphi _{2}\comp \chi \comp
\varphi _{1}^{-1})_{\mid W}$, for $k=1,...,n_{2}.$ Therefore,%
\begin{equation*}
(z_{1}\comp \varphi _{2}\comp \chi \comp \varphi
_{1}^{-1},...,z_{n_{2}}\comp \varphi _{2}\comp \chi \comp \varphi
_{1}^{-1})_{\mid W}=(f_{1},...,f_{n_{2}})_{\mid W}:W\rightarrow \mathbb{R}^{n_{2}},
\end{equation*}%
and $(f_{1},...,f_{n_{2}}):\mathbb{R}^{n_{1}}\rightarrow \mathbb{R}^{n_{2}}$
is an extension of $(\varphi _{2}\comp \chi \comp \varphi _{1}^{-1})_{\mid
W}:W\rightarrow V_{\varphi _{2}}$. The preceding argument can be repeated for every $x\in S_{1}$. Therefore, the
map $\chi :S_{1}\rightarrow S_{2}$ is smooth as a map of Aronszajn
subcartesian spaces.\hfill $\square $ \medskip 

To prove the converse to proposition 3, we need to
determine a differential structure for an Aronszajn subcartesian space. Let 
$(S,\mathfrak{A})$ be a subcartesian space of Aronszajn. The atlas $%
\mathfrak{A}$ on $S$ consists of charts $\varphi :U_{\varphi }\rightarrow
V_{\varphi }$, where $U_{\varphi }$ is an open subset of $S,$ $V_{\varphi }$
is a subset of $\mathbb{R}^{n_{\varphi }}$, for some $n_{\varphi }\in 
\mathbb{N}$, and $\varphi $ is a diffeomorphism of $U_{\varphi }$ onto $%
V_{\varphi }$. We assume that $\{U_{\varphi }\mid \varphi \in \mathfrak{A%
}\}$ is an open cover of $S$. Let $\mathcal{F}$ be the following family of functions
on the Aronszajn subcartesian space $(S,\mathfrak{A})$. 
A function $f:S\rightarrow \mathbb{R}$ is in $\mathcal{F}$, if for every $x\in S$ there exists 
$\varphi \in \mathfrak{A,}$ such that if $x\in U_{\varphi }$, then there exists an open neighbourhood 
$W_{\varphi (x)}$ of $\varphi (x)$ in $\mathbb{R}^{n_{\varphi }},$ and a function $F_{\varphi (x)}\in 
C^{\infty }\left( \mathbb{R}^{n_{\varphi }}\right) $ satisfying the
condition
\begin{equation}
f_{\mid \varphi ^{-1}(V_{\varphi }\cap W_{\varphi (x)})}=(F_{\varphi (x)}\comp \varphi )_{\mid
\varphi ^{-1}(V_{\varphi }\cap W_{\varphi (x)}).\smallskip }  
\label{F(S)}
\end{equation}

\noindent \textbf{Proposition 4} The family of sets 
\begin{equation*}
\mathfrak{s}=\{f^{-1}(I)\mid f\in \mathcal{F}\text{ and }I\text{ open
interval in }\mathbb{R}\}
\end{equation*}%
is a subbasis for the topology of the Aronszajn subcartesian space $(S,\mathfrak{A})$.\medskip 

\noindent \textbf{Proof.}
Every $U_{\varphi }\in \mathfrak{A}$ is an open set in $S$. It is homeomorphic to a
subset $V_{\varphi }\subseteq \mathbb{R}^{n_{\varphi }}$ with topology
induced by by the inclusion map $V_{\varphi }\hookrightarrow \mathbb{R}^{n_{\varphi }}$.
Let $C^{\infty }(V_{\varphi })$ be the
differential structure of $V_{\varphi }$ generated by its inclusion into $\mathbb{R}^{n_{\varphi }}$. In other words, 
$C^{\infty }(V_{\varphi })$ is generated by the family $\mathcal{F}_{V_{\varphi }}$ of functions 
$h:V_{\varphi }\rightarrow \mathbb{R}$ such that for every $\varphi (x)\in
V_{\varphi }$, there exists an open neighbourhood $W_{\varphi (x)}$ of $\varphi (x)$
in $\mathbb{R}^{n_{\varphi }}$ and a function $F_{\varphi (x)}\in C^{\infty}\left( \mathbb{R}^{n_{\varphi }}\right) $, which satisfies the condition%
\begin{equation*}
h_{\mid V_{\varphi }\cap W_{\varphi (x)}}=F_{\varphi (x) \mid V_{\varphi }\cap W_{\varphi (x)}}.
\end{equation*}
The topology of $(V_{\varphi },C^{\infty }(V_{\varphi }))$ has a subbasis 
\begin{equation*}
\mathfrak{s}_{V_{\varphi }}=\{h^{-1}(I)\mathfrak{\mid }h\in \mathcal{F}%
_{V_{\varphi }}\text{ and }I\text{ is an open interval in }\mathbb{R}\}.
\end{equation*}%
Since $\varphi :U_{\varphi }\rightarrow V_{\varphi }$ is a homeomorphism, 
it follows that $\varphi ^{-1}(\mathfrak{s}_{V_{\varphi }})$ is a subbasis for the
topology of $U_{\varphi }$. But 
\begin{eqnarray*}
\varphi ^{-1}(\mathfrak{s}_{V_{\varphi }}) &=&\{\varphi ^{-1}(h^{-1}(I))%
\mathfrak{\mid }h\in \mathcal{F}_{V_{\varphi }}\text{ and }I\text{ is an
open interval in }\mathbb{R}\} \\
&=&\{(h\comp \varphi )^{-1}(I)\mid h\in \mathcal{F}_{V_{\varphi }}\text{ and 
}I\text{ is an open interval in }\mathbb{R}\} \\
&\subseteq &\{f^{-1}(I)\mid f\in \mathcal{F}\text{ and }I\text{ is an open
interval of }R\}=\mathfrak{s.}
\end{eqnarray*}%
Hence 
\begin{equation}
\bigcup_{\varphi \in \mathfrak{A}}\varphi ^{-1}(\mathfrak{s}%
_{V_{\varphi }})\subseteq \mathfrak{s.}  \label{subbasis inclusion}
\end{equation}
Since $S=\cup _{\varphi \in \mathfrak{A}}U_{\varphi }$ and $\varphi ^{-1}(%
\mathfrak{s}_{V_{\varphi }})$ is a subbasis for the topology of $U_{\varphi} $, it follows that 
$\cup _{\varphi \in \mathfrak{A}}\varphi ^{-1}(%
\mathfrak{s}_{V_{\varphi }})$ is a subbasis for the topology of $S.$ The
inclusion (\ref{subbasis inclusion}) ensures that $\mathfrak{s}$ is a
subbasis for the topology of $S$.\hfill $\square $ \medskip 

There is a differential structure $C^{\infty }(S)$ on a set $S$, which is generated by a
family $\mathcal{F}$ of functions on $S$, see \cite{sniatycki}. Applying this construction to a
subcartesian space $(S,\mathfrak{A})$ of Aronszajn, and taking the family $\mathcal{F}$ defined by equation (\ref{F(S)}), we get the \emph{differential structure} $C^{\infty }(S)$ on $S$ \emph{determined by the atlas} 
$\mathfrak{A}$ on $S$. Proposition 4 ensures that the subcartesian space $(S,\mathfrak{A})$ and 
the corresponding differential space $(S,C^{\infty }(S))$ have the same topology. \medskip 

\noindent \textbf{Proposition 5} Let $(S_{1},\mathfrak{A}_{1})$ and $(S_{2},\mathfrak{A}_{2})$ be an Aronszajn subcartesian spaces, and let $(S_{1},C^{\infty }(S_{1}))$ and $(S_{2},C^{\infty }(S_{2}))$ be
the corresponding subcartesian differential spaces. If $\chi :(S_{1},%
\mathfrak{A}_{1})\rightarrow (S_{2},\mathfrak{A}_{2})$ is a smooth map of
subcartesian spaces of Aronszajn, then $\chi :(S_{1},C^{\infty}(S_{1}))\rightarrow (S_{2},C^{\infty }(S_{2}))$ is smooth map of subcartesian differential spaces. \medskip 

\noindent \textbf{Proof.}
By definition of Aronszajn subcartesian space, the assumption that the map 
$\chi :(S_{1},C^{\infty}(S_{1}))\rightarrow (S_{2},C^{\infty }(S_{2}))$ is smooth, means that for
every $x_{1}\in S_{1}$, there exist $\varphi _{1}:U_{\varphi
_{1}}\rightarrow V_{\varphi _{1}}$ in $\mathfrak{A}_{1}$ and $\varphi
_{2}:U_{\varphi _{2}}\rightarrow V_{\varphi _{2}}$ in $\mathfrak{A}_{2}$,
such that $x_{1}\in U_{\varphi _{1}}$, $x_{2}=\chi (x_{1})\in U_{\varphi
_{2}}$ and $\varphi _{2}\comp \chi \comp \varphi _{1}^{-1}:V_{\varphi
_{1}}\rightarrow V_{\varphi _{2}}$ extends to a $C^{\infty }$-map of a
neighbourhood of $\varphi _{1}(x_{1})\in \mathbb{R}^{n_{1}}$ to a
neighbourhood of $\varphi _{2}(x_{2})\in \mathbb{R}^{n_{2}}$, where $%
n_{1}=n_{\varphi _{1}}$ and $n_{2}=n_{\varphi _{2}},$ as in the proof of
proposition 2.\smallskip

For $i=1,2,$ let $\mathcal{F}_{i}$ denote the space of functions on $S_{i}$
determined by equation (\ref{F(S)}). It suffices to show that, for every $%
f_{2}\in \mathcal{F}_{2}$ the pull-back $f_{1}=\chi ^{\ast }f_{2}=f_{2}\comp
\chi \ $is in $\mathcal{F}_{1}$. In other words, we have to show that, for  
every $x_{1}\in S_{1}$, there exists $\varphi _{1}\in \mathfrak{A}_{1}$ such that if $x_{1}\in U_{\varphi _{1}}$ 
then there exists an open neighbourhood $W_{{\varphi }_1(x_{1})}$ of $\varphi _{1}(x_{1})$ in 
$\mathbb{R}^{n_{1}}$ and a function $F_{{\varphi }_1(x_{1})}\in C^{\infty }\left( \mathbb{R}^{n_{1}}\right) $ 
satisfying %
\begin{equation}
f_{1\mid \varphi _{1}^{-1}(V_{\varphi _{1}}\cap W_{{\varphi }_1(x_{1})})}=(F_{{\varphi }_1(x_{1})}\comp
\varphi _{1})_{\mid \varphi _{1}^{-1}(V_{\varphi _{1}}\cap W_{{\varphi }_1(x_{1})}).}  
\label{f1}
\end{equation}
On the other hand, $f_{1}=f_{2}\comp \chi $, where $f_{2}\in \mathcal{F}_{2}$. 
Therefore, for $x_{2}=\chi (x_{1})\in S_{2}$, there there exists $\varphi _{2}\in \mathfrak{A}_{2}$ such that 
if $x_{2}\in U_{\varphi _{2}}$, then there exists an open neighbourhood $W_{{\varphi }_2(x_{2})}$ of 
$\varphi _{2}(x_{2})$ in $\mathbb{R}^{n_{2}},$ and a function $F_{{\varphi }_2(x_{2})}\in C^{\infty}
\left( \mathbb{R}^{n_{2}}\right) $ satisfying %
\begin{equation}
f_{2\mid \varphi _{2}^{-1}(V_{\varphi _{2}}\cap W_{{\varphi }_2(x_{2})})}= 
(F_{{\varphi }_2(x_{2})}\comp \varphi _{2})_{\mid \varphi _{2}^{-1}(V_{\varphi _{2}}\cap W_{{\varphi }_2(x_{2})}).}
\label{f2}
\end{equation}%
By hypothesis $\chi $ is a smooth map of the Aronszajn subcartesian space $(S_1, {\mathfrak{A}}_1)$ 
into the Aronszajn subcartesian space $(S_2, {\mathfrak{A}}_2)$. So the map $\psi =
\varphi _{2}\comp \chi \comp \varphi _{1}^{-1}:V_{\varphi
_{1}}\rightarrow V_{\varphi _{2}}$ extends to a $C^{\infty }$-map $\widetilde{%
\psi}:\widetilde{W}_{1}\rightarrow \widetilde{W}_{2}$, where $\widetilde{W}_{1}$ is a
neighbourhood of $\varphi _{1}(x_{1})\in \mathbb{R}^{n_{1}}$, and $\widetilde{W} _{2}$ is a neighbourhood of 
$\varphi _{2}(x_{2})=\varphi _{2}(\chi (x_{1}))\in \mathbb{R}^{n_{2}}$. Without loss of generality, we may shrink 
$W_{{\varphi }_1(x_{1})},W_{{\varphi }_2(x_{2})},\widetilde{W}_{1}$ and $\widetilde{W}_{2}$ so that 
$W_{{\varphi }_1(x_{1})}=\widetilde{W}_{1}$ and $W_{{\varphi }_2(x_{2})}=\widetilde{W}_{2}$. Note, that the existence of appropriate $W_{{\varphi }_2(x_{2})},\widetilde{W}_{1}$ and $\widetilde{W}_{2}$ is guaranteed by
the hypotheses that $f_{2}$ and $\chi $ are smooth. The existence of an
appropriate $W_{{\varphi }_1(x_{1})}$ has to be established. Initially, we choose 
$W_{{\varphi }_1(x_{1}) }=\widetilde{W}_{1}=\widetilde{\psi}^{-1}(\widetilde{W}_{2})=\widetilde{\psi}%
^{-1}(W_{x_{2}}).$ Later, we shall have to shrink $W_{{\varphi }_1(x_{1}) }$ some more. \smallskip 

Having made these choices, we compute%
\begin{align}
(F_{{\varphi }_2(x_{2})}\comp \varphi _{2})_{\mid \varphi _{2}^{-1}(V_{\varphi _{2}}\cap
W_{{\varphi }_2(x_{2})})} &= F_{{\varphi }_2(x_{2}) \mid V_{\varphi _{2}}\cap W_{{\varphi }_2(x_{2})}}= 
F_{{\varphi }_2(x_{2})\mid \psi (V_{\varphi _{1}})\cap W_{{\varphi }_2(x_{2})}} \notag \\
&\hspace{-1in} =F_{{\varphi }_2(x_{2})\mid \widetilde{\psi}(V_{\varphi _{1}})\cap W_{{\varphi }_2(x_{2})}}  
=F_{{\varphi }_2(x_{2})\mid \widetilde{\psi}(V_{\varphi _{1}}\cap 
\widetilde{\psi}^{-1}(W_{{\varphi }_2(x_{2})}))} \notag \\
&\hspace{-1in} =(F_{{\varphi }_2(x_{2})}\comp \widetilde{\psi})_{\mid V_{\varphi _{1}}\cap 
\widetilde{\psi}^{-1}(W_{{\varphi }_2(x_{2})})}  
=(F_{{\varphi }_2(x_{2})}\comp \widetilde{\psi})_{\mid V_{\varphi _{1}}\cap
W_{{\varphi }_1(x_{1})}} \notag \\
&\hspace{-1in} = \big( ( F_{{\varphi }_2(x_{2})}\comp \widetilde{\psi}) \comp \varphi
_{1} \big)_{\mid \varphi _{1}^{-1}(V_{\varphi _{1}}\cap W_{{\varphi }_1(x_{1})}).}  
\label{f3}
\end{align}
Since $\psi =\varphi _{2}\comp \chi \comp \varphi _{1}^{-1}$ and $f_{1}=f_{2}\comp \chi ,$ 
\begin{align}
f_{1\mid \varphi _{1}^{-1}(V_{\varphi _{1}}\cap W_{{\varphi }_1(x_{1})})} &= 
(f_{2}\comp \chi )_{\mid \varphi _{1}^{-1}(V_{\varphi _{1}}\cap W_{{\varphi }_1(x_{1})})} \notag \\
& =(f_{2}\comp \varphi _{2}^{-1}\comp 
\psi \comp \varphi _{1})_{\mid \varphi _{1}^{-1}(V_{\varphi _{1}}\cap W_{{\varphi }_1(x_{1})})}  
\notag \\
& =(f_{2}\comp \varphi _{2}^{-1})_{\mid \psi (V_{\varphi _{1}}\cap W_{{\varphi }_1(x_{1})})}  
 =(f_{2}\comp \varphi _{2}^{-1})_{\mid (V_{\varphi _{2}}\cap W_{{\varphi }_2(x_{2})})} \notag \\ 
& =f_{2\mid \varphi _{2}^{-1}(V_{\varphi _{2}}\cap W_{{\varphi }_2(x_{2})})}.  
\label{f4}
\end{align}%
Taking into account equations (\ref{f1}) through (\ref{f4}), we get%
\begin{align}
f_{1\mid \varphi _{1}^{-1}(V_{\varphi _{1}}\cap W_{{\varphi }_1(x_{1}) })} & 
=f_{2\mid \varphi _{2}^{-1}(V_{\varphi _{2}}\cap W_{{\varphi }_2(x_{2})})}= 
(F_{{\varphi }_2(x_{2})}\comp \varphi _{2})_{\mid \varphi _{2}^{-1}(V_{\varphi _{2}}\cap W_{{\varphi}_2(x_{2})})} \notag \\
& =\big( (F_{{\varphi }_2(x_{2})}\comp \widetilde{\psi}) \comp 
\varphi _{1} \big)_{\mid \varphi _{1}^{-1}(V_{\varphi _{1}}\cap W_{{\varphi }_1(x_{1})})}.  
\label{f5}
\end{align}
Note that $F_{{\varphi }_2(x_{2})}\comp \widetilde{\psi}$ is a smooth function
of the open subset $\widetilde{\psi}^{-1}(W_{{\varphi }_2(x_{2}}))$ of $\mathbb{R}^{n_{1}}$.
However, it need not extend to a smooth function on $\mathbb{R}^{n_{1}}$.
But we can always find an open subset $W_{{\varphi }_1(x_{1})}$ of 
$\widetilde{\psi} ^{-1}(W_{{\varphi }_2(x_{2})})$ containing $x_{1}$ and $F_{{\varphi }_1(x_{1})} 
\in C^{\infty}(\mathbb{R}^{n_{1}})$ such that 
$F_{{\varphi }_1(x_{1})\mid W_{{\varphi }_1(x_{1})}}= (F_{{\varphi }_2(x_{2})}\comp 
\widetilde{\psi}) _{\mid W_{{\varphi}_1(x_{1})}}$. Equation (\ref{f5})
ensures that equation (\ref{f1}) is satisfied with the choices made here.
This argument works for every $f_{2}\in \mathcal{F}_{2}$. Hence 
$\chi :(S_{1},C^{\infty }(S_{1}))\rightarrow (S_{2},C^{\infty}(S_{2}))$ is smooth map of differential spaces. 
\hfill $\square $ \medskip 

In proving proposition 2 we have shown that to a given a subcartesian differential space $S$ with
differential structure $C^{\infty }(S)$, we can associate a
subcartesian space of Aronszajn with an atlas $\mathfrak{A}$ consisting of
diffeomorphisms $\varphi :U_{\varphi }\rightarrow V_{\varphi }\subseteq 
\mathbb{R}^{n_{\varphi }}$, where $U_{\varphi }$ is an open subset of $S$
and $V_{\varphi }$ is an arbitrary differential subspace of $\mathbb{R}^{n_{\varphi }}$. This construction gives a morphism $(S,C^{\infty}(S))\longrightarrow (S,\mathfrak{A})$ from the category of subcartesian
differential spaces to the category of subcartesian spaces of Aronszajn.
Conversely, given a subcartesian space of Aronszajn $(S,\mathfrak{A})$, in proving proposition 5 we
have constructed a morphism $(S,\mathfrak{A})\longrightarrow (S,C^{\infty }(S))$ from 
the category of subcartesian spaces of Aronszajn to the
category of subcartesian differential spaces. \medskip 

\noindent \textbf{Theorem 6} The morphisms 
\begin{displaymath}
(S,C^{\infty}(S))\longrightarrow (S,\mathfrak{A}) \, \, \, 
\mathrm{and} \, \, \,  (S,\mathfrak{A})\longrightarrow
(S,C^{\infty }(S)),
\end{displaymath} 
defined in the preceding paragraph, are inverses of each other. \medskip 

\noindent \textbf{Proof.}
Let $S$ be a subcartesian differential space with differential structure $%
C^{\infty }(S)$, and let $\mathfrak{A}$ be the atlas in $S$
consisting of diffeomorphisms 
\begin{displaymath}
\varphi :(U_{\varphi },C^{\infty}(U_{\varphi }))
\rightarrow (V_{\varphi },C^{\infty }(V_{\varphi })), 
\end{displaymath}
where $U_{\varphi }$, is an open subset of $S,$ $C^{\infty}(U_{\varphi })$ generated by the 
inclusion map $\iota :U_{\varphi }\hookrightarrow S$ and $(V_{\varphi },C^{\infty }(V_{\varphi }))$
is an arbitrary differential subspace of $\mathbb{R}^{n_{\varphi }}$. If $%
f\in C^{\infty }(S)$, then $f_{\mid U_{\varphi }}\in C^{\infty }(U_{\varphi })$ and 
$f_{\mid U_{\varphi }}\comp \varphi ^{-1}\in 
C^{\infty }(V_{\varphi })$. This means that, for every $x \in V_{\varphi }\subseteq $ 
$\mathbb{R}^{n_{\varphi }}$ there
exists an open neighbourhood $W_{x}$ of $x$ in $\mathbb{R}^{n_{\varphi }}$ and a smooth function 
$F_{x} \in C^{\infty }(\mathbb{R}^{n_{\varphi }})$ such that 
\begin{equation}
\left( f_{\mid U_{\varphi }}\comp \varphi ^{-1}\right) _{\mid V_{\varphi }\cap W_{x}} = 
F_{x}\mid {V_{\varphi }\cap W_{x}.}\smallskip  \label{F(S)a}
\end{equation}%
Conversely, a function $f:S\rightarrow \mathbb{R}$ such that, for every $%
\varphi \in \mathfrak{A}$, $f_{\mid U_{\varphi }}\comp \varphi ^{-1}\in 
C^{\infty }(V_{\varphi })$, then $f_{\mid U_{\varphi }}\in 
C^{\infty }(U_{\varphi })$. Since the differential structure $C^{\infty }(U_{\varphi })$ is generated by 
restrictions to $U_{\varphi }$ of functions in $C^{\infty }(S)$, the definition of 
differential structure ensures that $f\in C^{\infty }(S)$, see \cite{sniatycki}. \smallskip  

Let $\mathcal{F}$ be family of functions $f$ on $S$ such that,
for every $x\in S$, there exists $\varphi \in \mathfrak{A}$ with $U_{\varphi }$ an open subset of $S$ containing $x$, and there exists an open neighbourhood $W_{\varphi (x)}$
of $\varphi (x)$ in $\mathbb{R}^{n_{\varphi }},$ and a function $F_{\varphi (x)}\in 
C^{\infty }\left( \mathbb{R}^{n_{\varphi }}\right) $, which satisfies 
\begin{equation}
f_{\mid \varphi ^{-1}(V_{\varphi }\cap W_{\varphi (x) })}=(F_{\varphi (x)}\comp 
\varphi )_{\mid \varphi ^{-1}(V_{\varphi }\cap W_{\varphi (x)}).\smallskip } 
\label{F(S)b}
\end{equation}
The argument in the preceding paragraph shows that $\mathcal{F\subseteq C}^{\infty }(S).$ 
Proposition 4 ensures that the topology of the differential structure on 
$S$ generated by $\mathcal{F}$ coincides with the topopolgy of the original
differential structure $C^{\infty }(S)$. We need only show that the
differential structure generated by $\mathcal{F}$ coincides with original
differential structure $C^{\infty }(S)$. Since $\mathcal{F} \subseteq C^{\infty }(S)$, it follows that the differential structure generated by $\mathcal{F}$ is contained in $C^{\infty }(S)$. Conversely, let $f\in C^{\infty }(S)$. For every $\varphi \in \mathfrak{A}$ and $y\in V_{\varphi }$ it satisfies equation (%
\ref{F(S)a}). If $x=\varphi ^{-1}(y)$ and $F_{x}=F_y \comp \varphi ^{-1}$, then $f$ satisfies equation (\ref{F(S)b}).
So for every $\varphi \in \mathfrak{A}$, $f_{\mid U_{\varphi }}$ coincides with
the restriction to $U_{\varphi }$ of a function in the differential
structure generated by $\mathcal{F}$. This implies that $f_{\mid U_{\varphi }}$ lies in 
the differential structure generated by $\mathcal{F}$. Hence 
the differential structure on $S$ generated by $\mathcal{F}$ is equal to the
original differential structure $C^{\infty }(S)$. \hfill $\square $ \medskip

In the construction of the morphisms 
$(S,C^{\infty}(S))\rightarrow (S,\mathfrak{A})$ and $(S,\mathfrak{A})\rightarrow (S,C^{\infty }(S))$  
we used only intrinsic constructions. Therefore these morphisms are natural transformations. In particular, all
geometric structures and relations obtained in the category of subcartesian
spaces of Aronszajn may be translated to the category of subcartesian differential spaces.


\begin{thebibliography}{9}
\bibitem{aronszajn} N. Aronszajn, Subcartesian and subRiemannian spaces, \textit{Notices 
American Mathematical Society} \textbf{14} (1967) 111.

\bibitem{sikorski1} R. Sikorski, Abstract covariant derivative, \textit{%
Colloq. Math. }\textbf{18} (1967) 252-272.

\bibitem{sikorski72} R. Sikorski, \textit{Wst\k{e}p do geometrii r\'{o}\.{z}niczkowej} [Introduction to differential geometry], Pa\'{n}stwowe Wydawnictwo Naukowe, Warsaw, 1972. 

\bibitem{sniatycki} J. \'{S}niatycki, \textit{Differential geometry of singular spaces and 
reduction of symmetry}, Cambridge University Press, Cambridge, UK, 2013.

\bibitem{walczak}P. Walczak, A theorem on diffeomorphisms in the category of differential spaces. 
\textit{Bull. Acad. Polon. Sci. S\'{e}r. Sci. Math. Astronom. Phys.} \textbf{21} (1973) 325--329.


\end{thebibliography}
\end{document}